\documentclass[12pt,dvipdfm]{article}
\usepackage{amsfonts}
\usepackage{graphicx}
\usepackage{amsfonts, amsmath, amssymb, latexsym, epsfig}
\usepackage{mathrsfs}
\newcommand{\blst}{\begin{trivlist}}
\newcommand{\elst}{\end{trivlist}}

\newtheorem{thm}{Theorem}[section]
\newtheorem{prop}[thm]{Proposition}
\newtheorem{cor}[thm]{Corollary}
\newtheorem{lem}[thm]{Lemma}
\newtheorem{conj}[thm]{Conjecture}
\newtheorem{exa}[thm]{Example}
\newtheorem{defn}[thm]{Definition}

\newcommand{\ben}{\begin{enumerate}}
\newcommand{\een}{\end{enumerate}}
\newcommand{\ble}{\begin{lem}}
\newcommand{\ele}{\end{lem}}
\newcommand{\bth}{\begin{thm}}
\renewcommand{\eth}{\end{thm}}
\newcommand{\bpr}{\begin{prop}}
\newcommand{\epr}{\end{prop}}
\newcommand{\bco}{\begin{cor}}
\newcommand{\eco}{\end{cor}}
\newcommand{\bcon}{\begin{conj}}
\newcommand{\econ}{\end{conj}}
\newcommand{\bde}{\begin{defn}}
\newcommand{\ede}{\end{defn}}
\newcommand{\bex}{\begin{exa}}
\newcommand{\eex}{\end{exa}}
\newcommand{\barr}{\begin{array}}
\newcommand{\earr}{\end{array}}
\newcommand{\btab}{\begin{tabular}}
\newcommand{\etab}{\end{tabular}}
\newcommand{\beq}{\begin{equation}}
\newcommand{\eeq}{\end{equation}}

\newcommand{\bea}{\begin{eqnarray*}}
\newcommand{\eea}{\end{eqnarray*}}
\newcommand{\beaa}{\begin{eqnarray}}
\newcommand{\eeaa}{\end{eqnarray}}
\newcommand{\bce}{\begin{center}}
\newcommand{\ece}{\end{center}}
\newcommand{\bpi}{\begin{picture}}
\newcommand{\epi}{\end{picture}}
\newcommand{\bfi}{\begin{figure} \begin{center}}
\newcommand{\efi}{\end{center} \end{figure}}

\newcommand{\bsl}{\begin{slide}{}}
\newcommand{\esl}{\end{slide}}

{}
{}
{}
{}
{}
{}
{}
{}
{}
{}
{}
{}
{}
{}
{}
{}
{}
{}
{}
{}
{}
{}
{}
{}
{}
{}
{}
{}
\newenvironment{proof}{
\par
\noindent {\bf Proof.}\rm}{\mbox{}\hfill\rule{0.5em}{0.809em}\par}
\begin{document}
\title{Enumerations of Permutations by Circular Descent Sets}
\author{ Hungyung Chang$^{a,}$\thanks{Email address of the corresponding author: giannic@math.sinica.edu.tw}
 \and  Jun Ma$^b$\\
 \and  Yeong-Nan Yeh$^{c,}$\thanks{Partially supported by NSC 96-2115-M-006-012}}
\maketitle \vspace*{-1.2cm}
\date{}
\bce \footnotesize $^{a,b,c}$ Institute of Mathematics, Academia
Sinica, Taipei, Taiwan. \ece \vspace*{-0.3cm} \thispagestyle{empty}
\begin{abstract}
The circular descent of a permutation $\sigma$ is a set
$\{\sigma(i)\mid \sigma(i)>\sigma(i+1)\}$. In this paper, we focus
on the enumerations of permutations by the circular descent set. Let
$cdes_n(S)$ be the number of permutations of length $n$ which have
the circular descent set $S$. We derive the explicit formula for
$cdes_n(S)$. We describe a class of generating binary trees $T_k $
with weights. We find that the number of permutations in the set
$CDES_n(S)$ corresponds to the weights of $T_k$. As a application of
the main results in this paper, we also give the enumeration of
permutation tableaux according to their shape.
\end{abstract}
\newpage
\section{Introduction}
Throughout this paper, let $[m,n]:=\{m,m+1,\cdots,n\}$, $[n]:=[1,n]$
and $\mathfrak{S}_n$ be the set of all the permutations in the set
$[n]$. We will write permutations of $\mathfrak{S}_n$ in the form
$\sigma=(\sigma(1)\sigma(2)\cdots\sigma(n))$. We say that a
permutation $\sigma$ has a {\it circular descent} of value
$\sigma(i)$ if $\sigma(i)>\sigma(i+1)$. The {\it circular descent
set} of a permutation $\sigma$, denoted $CDES(\sigma)$, is the set
$\{\sigma(i)\mid \sigma(i)>\sigma(i+1)\}.$ For example, the circular
descent set of $\sigma=(48632517)$ is $\{8,6,5,3\}$. For any
$S\subseteq [n]$, we define a set $CDES_n(S)$ as
$CDES_n(S)=\{\sigma\in\mathfrak{S}_n\mid CDES(\sigma)=S\}$ and use
$cdes_n(S)$ to denote the number of the permutations in the set
$CDES_n(S)$, i.e., $cdes_n(S)=|CDES_n(S)|$. First, for any $n\geq 2$
and $S\subseteq [n]$, we prove that the necessary and sufficient
condition for $CDES_n(S)\neq \emptyset$ is $1\notin S$. Furthermore,
 suppose that
$S=\{s_1,s_2,\ldots,s_k\}$ with $s_1>s_2>\ldots
>s_k>s_{k} > 1$. Let $d_i(S)=s_i-s_{i+1}$ for all $i\in [k-1]$ and $d_k(S)= s_k-1$.
The main results in this paper are as follows:
 \begin{eqnarray*}cdes_n(S) =
\displaystyle{\sum_{\substack{ x_1, \cdots, x_{|S|} \in \{0,1\}}}
(-1)^{|S|-\sum\limits_{j=1}^{|S|}x_j} \prod_{i=1}^{|S|}
\left(1+\sum\limits_{j=1}^ix_i\right)^{d_i(S)} }\end{eqnarray*}
for each positive integer $n$ and $S\subseteq [2,n]$.

The concept of generating trees has been introduced in the
literature by Chung, Graham, Hoggat and Kleiman in \cite{FRKC} to
examine Baxter permutations. There are closed relations between the
generating trees and many combinatorial models. Also we describe a
class of generating binary trees $T_k $ with weights. We find that
the number of permutations in the set $CDES_n(S)$ corresponds to the
weights of $T_k $.

A. Postnikov \cite{post2006} and L. Williams \cite{williams}
mentioned the conceptions of permutation tableau when they studied
the enumeration of the totally positive Grassmannian cells.
Surprisingly, there are the closed connections between permutation
tableaux and a statistical physics model called the Partially
Asymmetric Exclusion Process \cite{corteel1,corteel2,corteel3}.

Following \cite{ESLK}, we give the definition of the permutation
tableau as follows. Recall that a partition
$\lambda=(\lambda_1,\cdots,\lambda_k)$ is a weakly decreasing
sequence of non-negative integers. For a partition $\lambda$, which
$\Sigma \lambda_i = m$, the Young diagram $Y_\lambda$ of Shape
$\lambda$ is a left-justified diagram of $m$ boxes, with $\lambda_i$
boxes in the $i$th row.

Define a permutation tableau ${\cal T}^k_n$ to be a partition
$\lambda$ such that $Y_{\lambda}$ is contained in a $k\times (n-k)$
rectangle, together with filling of the boxes of $Y_\lambda$ with
$0's$ and $1's$ such that the following properties hold:

(1) Each column of the rectangle contains at least one $1$.

(2) There is no $0$ which has a $1$ above it in the same column and
a $1$ to its left in the same rows.

The partition $\lambda$ is called the shape of the permutation
tableau ${\cal T}^k_n$. Using the bijection $I$ in \cite{SCPN}, as a
application of the main results, we give the enumeration of
permutation tableaux according to their shape.

A weak excedance of a permutation $\sigma$ is a value $\sigma(i)$
such that $\sigma(i)\geq i$. In this situation, $i$ ( resp.
$\sigma(i)$) are called a weak excedance bottom (resp. top) of
$\sigma$. Obviously, Non-weak excedance bottom and Non-weak
excedance top can be defined in terms of $i$ and $\sigma(i)$ such
that $\sigma(i)<i$. Furthermore, for any $\sigma\in\mathfrak{S}_n$,
the non-weak excedance bottom set of a permutation $\sigma$, denoted
by $NWEXB(\sigma)$, as $NWEXB(\sigma)=\{i\mid \sigma(i)<i\}$. For
any $S\subseteq [n]$, define a set $NWEXB_n(S)$ as
$NWEXB_n(S)=\{\sigma\in\mathfrak{S}_n\mid NWEXB(\sigma)=S\}$ and use
$nwexb_n(S)$ to denote the number of the permutations in the set
$NWEXB_n(S)$, i.e., $nwexb_n(S)=|NWEXB_n(S)|$.

 We must indicate that
there are closed relation between $nwexb_n(S)$ and the Genocchi
numbers. The study of Genocchi numbers and their combinatorial
interpretations has received much attention
\cite{Dom,DD1,DD2,DDF,DDR,DDV,GAN,Rio73}. Particularly, Michael
Domartzki \cite{Dom} give the definition of $k$-th generalized
Genocchi numbers $\{G^{(k)}_{2n}\}_{n\geq 1}$ as follows:

 Let $A^{(k)}_{n+1}(X)$ be the following Gandhi polynomial
in $X$:
\begin{eqnarray*}A^{(k)}_{n+1}(X)&=&X^kA^{(k)}_{n}(X+1)-(X-1)^kA^{(k)}_{n}(X)~~\forall n\geq 0\\
A^{(k)}_{0}(X)&=&1.\end{eqnarray*}Define $k$-th generalized Genocchi
numbers $\{G^{(k)}_{2n}\}_{n\geq 1}$ by
$G^{(k)}_{2n}=A^{(k)}_{n-1}(1)$. In \cite{DD2} , Dumont has given
several interpretations of the usual Genocchi numbers $G^{(2)}_{2n}$
in terms of permutations. Michael Domartzki \cite{Dom} conclude that
the number of permutations $\sigma$ of $[kn]$ such that
$\sigma(i)\geq i$ iff $\sigma(i)\equiv 0(mod~ k)$ is equal to
$G_{2n+2}^{(k)}$. As another direct application of the main results
in this paper, we give the enumeration of permutations by the
non-weak excedance bottom set. This generalizes the results in
\cite{Dom}.

We organize this paper as follows. In Section 2, we study the
recursions and the generating functions for the sequence
$cdes_n(S)$. In Section 3, we give the proof of the main theorem.
and consider a class of the generating binary trees. In Section 4,
we give some applications of the main theorem of this paper. We
enumerate permutation tableaux according to their shape. We also
counting the number of the permutations by the non-weak excedance
bottom set.

\section{The Recursions and The Generating Functions  }
In this section, we will derive some recursions  and consider the
generating functions for the sequence $cdes_n(S)$. We start with
some definitions. For any $\sigma\in\mathfrak{S}_n$, let
$\tau=\sigma(i_1)\sigma(i_2)\ldots\sigma(i_k)$ be a subsequence of
$\sigma$ and ``red" the increasing bijection of
$\{\sigma(i_1),\sigma(i_2),\ldots,\sigma(i_k)\}$ onto $[k]$. Define
the reduction of the subsequence $\tau$, denote $red(\tau)$, as
$red(\sigma(i_1))red(\sigma(i_2))\ldots red(\sigma(i_k))$. For any
$S\subseteq [2,n]$, define a set $\delta(S)$ as $\delta(S)=\{s-1\mid
s\in S\}$, clearly, $\delta(S)\subseteq [n-1]$. First, we easily
obtain the following lemma.

\begin{lem}
\label{d1} Let $n\geq 1$ and $S\subseteq [n]$. Then $CDES_n(S) \neq
\emptyset$ if and only if  $1\notin S$.
\end{lem}
\begin{proof}
Note that $1$ is impossible to be a circular descent. For any
$S\subseteq [2,n]$, suppose that $S=\{i_1,\ldots,i_k\}$, where
$i_1<i_2<\ldots <i_k$. Let $T=[n]\setminus
S=\{1,j_1,\ldots,j_{n-k-1}\}$, where $1<j_1< \ldots <j_{n-k-1}$.
Then $(i_k,i_{k-1},\ldots ,i_1,1,j_1,\ldots,j_{n-k-1})\in
CDES_n(S)$.
\end{proof}

\begin{lem}
\label{DesR}Suppose $n$ is a positive integer with $n\geq 2$ and
$S\subseteq [2,n]$. Let $i=\min(S)$. Then $$ cdes_n(S) = cdes_n(S
\cup \{ i-1 \} \backslash \{ i \})+ cdes_{n-1}(\delta(S)) +
cdes_{n-1}(\delta(S) \backslash \{ i-1 \}).$$
\end{lem}
\begin{proof}
Let $T_1$ be a set of the permutations $\sigma$ in the set
$CDES_n(S)$ such that $i$ and $i-1$ are not adjacent in $\sigma$,
i.e., $T_1 = \{ \sigma \in CDES_n(S) \mid
|\sigma^{-1}(i)-\sigma^{-1}(i-1)|\geq 2 \}$. Let $T_2$ (resp. $T_3$)
be a set of the permutations $\sigma$ in the set $CDES_n(S)$ such
that the position of $i-1$ (resp. $i$) is exactly the left of $i$
(resp. $i-1$), i.e., $T_2 = \{ \sigma \in CDES_n(S) \mid
\sigma^{-1}(i)=\sigma^{-1}(i-1)+1 \}$ (resp. $T_3 = \{ \sigma \in
CDES_n(S) \mid \sigma^{-1}(i-1)=\sigma^{-1}(i)+1 \}$). Furthermore,
let $S_1 = S \cup \{ i-1 \} \backslash \{ i \}$ and $S_2 =
\delta(S)$ and $S_3 = \delta(S) \backslash \{ i-1 \}$. To obtain the
results of the lemma, it is sufficient to give bijections from $T_1$
to $CDES_n(S_1)$, $T_2$ to $CDES_{n-1}( S_2)$ and $T_3$ to
$CDES_{n-1}(S_3)$, respectively. We discuss the following two cases.

{\it Case 1.} The bijection from $T_1$ to $CDES_n(S_1)$

For any $\sigma\in T_1$, we exchange the positions of $i$ and $i-1$
in $\sigma$ and denote the obtained permutation by $\tau.$ Suppose
that $\sigma(j)=i$ and $\sigma(k)=i-1$. Then $i-1\notin S$ implies
that either $\sigma(k+1)>i-1$ or $k=n$. It is easy to check
$\sigma(j+1)<i-1$ since $i$ and $i-1$ are not adjacent in $\sigma$.
Hence, $i-1$ is a circular descent in $\tau$. Similarly, we may
prove that $i$ is impossible to be a circular descent in $\tau$.
Obviously, this is a bijection.

{\it Case 2.} The bijection $\phi$ from $T_2$ to $CDES_n(S_2)$

First, we claim that $T_2=\emptyset$ if $i=2$. Otherwise, the letter
$2$ can't be a circular descent since the position of $1$  is exact
the left of $2$. Note that $2\in S$ implies $1\in S_2$. Lemma
\ref{d1} tells us that $CDES_{n-1}(S_2)=\emptyset$. When $i\geq 3$,
we delete the letter $i$ from $\sigma$ and denote the obtained
subsequence as $\tau$. Let $\phi(\sigma)=red(\tau)$. Suppose
$\sigma(j)=i$. Then $\sigma(j+1)<i-1$ since  the position of $i-1$
is exactly the left of $i$. Hence, $i-1$ is a circular descent of
$\phi(\sigma)$ and $\phi(\sigma)\in CDES_{n-1}(S_2)$. Conversely,
for any $\tau\in CDES_{n-1}(S_2)$, we add $1$ to each of the letters
larger than $i-1$ and insert the letter $i$ behind $i-1$. Denote the
obtained sequence as $\phi^{-1}(\tau)$. Suppose that $\tau(j)=i-1$.
Then $\tau(j+1)<i-1$ since $i-1$ is a circular descent of $\tau$.
This implies $\phi^{-1}(\tau)\in T_2$, hence, $\phi$ is a bijection
from $T_2$ to $CDES_n(S_2)$.

The bijection from $T_3$ to $CDES_n(S_3)$ is similar to Case $2$.
The proof is completed.
\end{proof}

\begin{cor}
\label{Con2} Suppose $n$ is a positive integer with $n\geq 3$ and
$S\subseteq [3,n]$. Then $cdes_n(S\cup\{2\})=
cdes_{n-1}(\delta(S))$.
\end{cor}
\begin{proof}
For any $\sigma \in CDES_n(S\cup\{2\})$, the position of the letter
$2$ have to be exactly the left of $ 1$ in $\sigma$. Let $T_3$ be
defined as that in Lemma \ref{DesR}. Then $CDES_n(S\cup\{2\})=T_3$.
Hence, $cdes_n(S\cup\{2\})=cdes_{n-1}( \delta(S\cup\{2\})\backslash
\{1\} ) = cdes_{n-1}(\delta(S))$.
\end{proof}

\begin{lem}\label{lemmajn}
Suppose $n$ is a positive integer with $n\geq 2$ and $S\subseteq
[2,n]$. Let $j=\max(S)$. Then \\
(1) $cdes_j(S) = cdes_n(S)$ for each integer $n \geq j$. \\
(2) $cdes_{n}( \emptyset ) = 1$ for each positive integer $n$.
\end{lem}
\begin{proof} For any $\sigma \in
CDES_n(S)$,  since $ j+1, \cdots, n $ are not circular descents of
$\sigma$,  we have $\sigma(k)=k$ for all $k \in [j+1,n]$. This
implies  $cdes_{n}(S)=cdes_{j}(S)$. Similarly, we have $cdes_{n}(
 \emptyset ) = 1$.
\end{proof}

\begin{lem}\label{ex2}Suppose $n$ is a positive integer with $n\geq 2$. Then $cdes_n(\{n\})=2^{n-1}-1$
\end{lem}
\begin{proof}
Let $A =\{ \sigma(t) : t < \sigma^{-1}(n) \}$ and $B = \{ \sigma(t)
: t > \sigma^{-1}(n) \}$. Then $B \neq \emptyset$ since $n$ is the
circular descent of $\sigma$.  Furthermore, the elements of $A$ and
$B$ form a increasing subsequence of $\sigma$ repectively since $n$
are the unique circular descent.  For each letter $j\neq n$, the
position of $j$ has two possibilities at the left or right of $n$.
Therefore, $cdes_n(\{n\})=2^{n-1}-1$ since $B\neq \emptyset$.
\end{proof}

In the next lemma, we will derive another recursion for the sequence
$cdes_n(S)$.
\begin{lem}
\label{Yeh} Suppose $n$ is a positive integer with $n\geq 3$ and
$S\subseteq [2,n-1]$. Then  $$ \displaystyle{
cdes_n(S\cup\{n\})=(n-1-|S|)cdes_{n-1}(S) + \sum_{i \in [n-1]
\backslash S } cdes_{n-1}( S \cup \{ i \}) }.$$
\end{lem}
\begin{proof}
Suppose $i = \max (S)$. It is clear that $i < n$.   For each $\sigma
\in \mathfrak{S}_{n-1}$, we may extend $\sigma$ to be a permutation
of $\mathfrak{S}_n$ by inserting $n$ into $\sigma$. Therefore, there
is at most one descent of $\sigma$ that is eliminated by inserting
$n$. It is easy to see that there are $n-1-|S|$ ways to insert $n$
such that the set of descent of the new permutation is $S$. It is
clear that we have put $i \rightarrow_* n$ in the new permutation if
we want to eliminate $i \in S$. So we have
$cdes_n(S\cup\{n\})=(n-1-|S|)cdes_{n-1}(S) + \sum\limits_{i \in
[n-1]/S } cdes_{n-1}( S \cup \{ i \})$.
\end{proof}

For any $S\in [2,n]$, suppose $S=\{i_1,i_2,\ldots,i_k\}$, let
$\mathbf{x}_{S}$ stand for the monomial $x_{i_1-1}x_{i_2-1}\cdots
x_{i_k-1}$; particularly, let $\mathbf{x}_{\emptyset}=1$. Given
$n\geq 2$, we define a generating function
$${g_n(x_1,x_2,\ldots,x_{n-1};y) = \sum\limits_{\sigma \in
\mathfrak{S}_{n}} \mathbf{x}_{CDES(\sigma)}}y^{|CDES(\sigma)|}.$$ We
also write $g_n(x_1,x_2,\ldots,x_{n-1};y)$ as $g_n$ for short.  We
call it {\em circular descent polynomial}.

\begin{cor}Let $n$ be a positive integer with $n\geq 2$ and $g_n=\sum\limits_{\sigma \in
\mathfrak{S}_{n}} \mathbf{x}_{CDES(\sigma)}y^{|CDES(\sigma)|}.$ Then
$g_n$ satisfies the following recursion:
\begin{eqnarray*}g_{n+1}=(1+nx_ny)g_n+x_n\sum\limits_{i=1}^{n-1}\frac{\partial
g_n}{\partial x_i}-x_ny^2\frac{\partial g_n}{\partial
y}\end{eqnarray*} for all $n\geq 2$ with initial condition
$g_2=1+x_1y$, where the notation $\frac{\partial g_n}{\partial y}$
denote partial differentiation of $g_n$ with respect to $y$.
\end{cor}
\begin{proof} Obviously, $g_2=1+x_1y$ and  ${\sum\limits_{\sigma \in
\mathfrak{S}_{n}} \mathbf{x}_{CDES(\sigma)}y^{|CDES(\sigma)|} =
\sum\limits_{S\subseteq [2,n]}cdes_{n}(S)\mathbf{x}_{S}}y^{|S|}.$
Hence,
\begin{eqnarray*}g_{n+1}
&=&\sum\limits_{S\subseteq [2,n+1]}cdes_{n+1}(S)\mathbf{x}_{S}y^{|S|}\\
&=&\sum\limits_{S\subseteq [2,n+1],n+1\in
S}cdes_{n+1}(S)\mathbf{x}_{S}y^{|S|}
+\sum\limits_{S\subseteq [2,n+1],n+1\notin S}cdes_{n+1}(S)\mathbf{x}_{S}y^{|S|}\\
&=&\sum\limits_{S\subseteq
[2,n]}\left[(n-|S|)cdes_n(S)+\sum\limits_{i\in[2,n]\setminus
S}cdes_n(S\cup\{i\})\right]\mathbf{x}_{S}x_ny^{|S|+1}
+g_n\\
&=&\sum\limits_{S\subseteq [2,n]}\sum\limits_{i\in[2,n]\setminus
S}cdes_n(S\cup\{i\})\mathbf{x}_{S}x_ny^{|S|+1}
-\sum\limits_{S\subseteq
[2,n]}|S|cdes_n(S)\mathbf{x}_{S}x_ny^{|S|+1} +(1+nx_ny)g_n\\.
\end{eqnarray*}
Note that \begin{eqnarray*}\frac{\partial g_n}{\partial
y}=\sum\limits_{S\subseteq
[2,n]}|S|cdes_n(S)\mathbf{x}_{S}y^{|S|-1}\end{eqnarray*} and
\begin{eqnarray*}&&\sum\limits_{S\subseteq [2,n]}\sum\limits_{i\in[2,n]\setminus
S}cdes_n(S\cup\{i\})\mathbf{x}_{S}x_ny^{|S|+1}\\
&=&\sum\limits_{S\subseteq [2,n],S\neq
\emptyset}cdes_n(S)x_ny^{|S|}\sum\limits_{i\in
S}\frac{\mathbf{x}_{S}}{x_{i-1}}\\
&=&x_n\sum\limits_{i=1}^{n-1}\frac{\partial g_n}{\partial
x_i}.\end{eqnarray*} Therefore,
\begin{eqnarray*}g_{n+1}=(1+nx_ny)g_n+x_n\sum\limits_{i=1}^{n-1}\frac{\partial g_n}{\partial
x_i}-x_ny^2\frac{\partial g_n}{\partial y}.\end{eqnarray*}
\end{proof}

By computer search, we obtain the values of $cdes_n(S)$ for all
$2\leq n\leq 5$ and $S\subseteq [2,n]$. In the following table, we
list the generating function $g_n$ for $2\leq n\leq 5$.
\begin{center}
\begin{tabular}{|l|}
\hline \multicolumn{1}{|c|}{The Table of Circular Descent Polynomial
$g_n$ ( $ 2\leq n \leq 5 $ )} \\ \hline
 $g_2 = 1 + x_1y$ \\
 $g_3= 1 + x_1y + 3x_2y+ x_1x_2y^2$ \\
 $g_4 = 1 + x_1y + 3 x_2y + x_1x_2y^2 + 7x_3y + 3x_1x_3y^2 + 7x_2x_3y^2 + x_1x_2x_3y^3$ \\
 $g_5 = 1 + x_1y + 3 x_2y + x_1x_2y^2 + 7x_36 + 3x_1x_3y^2 + 7x_2x_3y^2 + x_1x_2x_3y^3 + 15 x_4y + 7 x_1 x_4y^2$ \\
  $~~~~~+ 17 x_2x_4y^2 + 3x_1x_2x_4y^3 + 31 x_3x_4y^2 + 7 x_1x_3x_4y^3 + 15 x_2x_3x_4y^3 + x_1x_2x_3x_4y^4$
       \\ \hline
\end{tabular}
\end{center}
\vspace*{1cm}
\section{Counting permutations by circular descent sets}
In this section, we will give the enumerations of permutations by
circular descent sets.
\begin{thm}\label{main} Let $n$ be a positive integer and $S\subseteq [2,n]$. Then
 \begin{eqnarray*}cdes_n(S) =
\displaystyle{\sum_{\substack{ x_1, \cdots, x_{|S|} \in \{0,1\}}}
(-1)^{|S|-\sum\limits_{j=1}^{|S|}x_j} \prod_{i=1}^{|S|}
\left(1+\sum\limits_{j=1}^ix_i\right)^{d_i(S)} }.\end{eqnarray*}
\end{thm}
\begin{proof}
Let $k=|S|$. By induction on $k$. For convenience, we let
$$R_k(S;x_1,\ldots,x_k)=
\displaystyle{(-1)^{k-x_1-\cdots-x_{k}}\prod_{j=1}^{k}
(1+x_1+\cdots+x_j)^{d_j(S)}} .$$  By Lemma \ref{lemmajn}(1), we may
always suppose $n\in S$. When $k=1$, $$cdes_{n}(S) =
\displaystyle{\sum_{\substack{ x_1 \in \{0,1\}}}
(-1)^{1-x_1}(1+x_1)^{n-1} } = 2^{n-1} -1.$$ By Lemma \ref{ex2}, the
theorem holds when $k=1$. Now, assume that the theorem holds for
$k'\leq k$, i.e.,
$$cdes_{n}(S) =\displaystyle{\sum_{\substack{ x_1, \cdots, x_{k'} \in
\{0,1\}}} R_{k'}(S;x_1,\ldots,x_{k'}) }.$$ We consider the case with
$k'=k+1$. If $s_{k+1}=2$, then let $S'=S\backslash\{s_{k+1}\}$.
Lemma \ref{Con2} tells us that $cdes_n( S )=cdes_n(
S'\cup\{2\})=cdes_{n-1}(\delta(S'))$. Furthermore, we have
$d_j(\delta(S'))=d_j(S)$ for any $j\in[k]$ and $d_{k+1}(S)=1$.  Note
that $|\delta(S')|=k$. So,
\begin{eqnarray*}R_{k+1}(S;x_1,\ldots,x_k,1)&=&R_k(\delta(S');x_1,\ldots,x_k)\left(2+\sum\limits_{i=2}^{k}x_i\right)\\
R_{k+1}(S;x_1,\ldots,x_k,0)&=&-R_k(\delta(S');x_1,\ldots,x_k)\left(1+\sum\limits_{i=2}^{k}x_i\right)\end{eqnarray*}
Therefore,
$$R_{k+1}(S;x_1,\ldots,x_k,1)+R_{k+1}(S;x_1,\ldots,x_k,0)=R_k(\delta(S');x_1,\ldots,x_k).$$
By induction hypothesis, we have
$$cdes_{n}(\delta(S')) = \displaystyle{\sum_{\substack{ x_1, \cdots,
x_k \in \{0,1\}}} R_k(\delta(S');x_1,\ldots,x_k) }.$$ Hence,
\begin{eqnarray*}cdes_{n}(S) &=&cdes_{n}(\delta(S'))\\
&=&\displaystyle{\sum_{\substack{ x_1, \cdots, x_k \in \{0,1\}}}
R_k(\delta(S');x_1,\ldots,x_k) }\\
&=&\displaystyle{\sum_{\substack{ x_1, \cdots, x_k \in \{0,1\}}}
\left[R_{k+1}(S;x_1,\ldots,x_k,1)+R_{k+1}(S;x_1,\ldots,x_k,0)\right]}\\
 &=& \displaystyle{\sum_{\substack{ x_1, \cdots, x_{k+1} \in \{0,1\}}}
R_{k+1}(S;x_1,\ldots,x_{k+1}) }.\end{eqnarray*}So, the theorem holds
for $k'=k+1$ and $i=2$. Use induction on $i$ again. Suppose the
theorem holds for $k'= k+1$ and $ 2 \leq i' \leq i$. Suppose  $i' =
i+1$. Let the sets $S_j$ be defined as that in Lemma \ref{DesR} for
$j=1,2,3$. It is easy to obtain that $d_j(S_1)=
d_j(S_2)=d_j(S_3)=d_j(S)$ for any $j \in[k-1]$, $d_k(S_1) =
d_k(S)+1$ , $d_{k+1}(S_1) =d_{k+1}(S_2)= d_{k+1}(S)-1$, $d_k(
S_2)=d_k(S)$,  and $d_k(S_3)=d_k(S)+d_{k+1}(S)-1$. Hence, we have
\begin{eqnarray*}   R_{k+1}(S_1;x_1,\ldots,x_{k+1})
&=& R_{k+1}(S;x_1,\ldots,x_{k+1})
\left(1+\sum\limits_{j=1}^{k}x_j\right)\left(1+\sum\limits_{j=1}^{k+1}x_j\right)^{-1}\\
R_{k+1}(S_2;x_1,\ldots,x_{k+1}) &=& R_{k+1}(S;x_1,\ldots,x_{k+1})
\left(1+\sum\limits_{j=1}^{k+1}x_j\right)^{-1}\\
 R_{k}(S_3;x_1,\ldots,x_{k})
&=&
-\left(1+\sum\limits_{j=1}^{k}x_j\right)^{-1}R_{k+1}(S;x_1,\ldots,x_{k},0)
\end{eqnarray*}
Note that
\begin{eqnarray*}&&R_{k+1}(S;x_1,\ldots,x_{k+1})-R_{k+1}(S_1;x_1,\ldots,x_{k+1})-R_{k+1}(S_2;x_1,\ldots,x_{k+1})\\
&=&R_{k+1}(S_1;x_1,\ldots,x_{k+1})\left[1-\left(1+\sum\limits_{j=1}^{k}x_j\right)\left(1+\sum\limits_{j=1}^{k+1}x_j\right)^{-1}-\left(1+\sum\limits_{j=1}^{k+1}x_j\right)^{-1}\right]\\
&=&\frac{x_{k+1}-1}{1+x_1+\ldots
+x_{k+1}}R_{k+1}(S;x_1,\ldots,x_{k+1}).\end{eqnarray*} Since
$min(S_1)=min(S_2)=i$, $|S_3|=k$, by induction hypothesis, we have
\begin{eqnarray*}  cdes_n(S_1) &=&\displaystyle{\sum_{\substack{ x_1, \cdots, x_{k+1} \in
\{0,1\}}} R_{k+1}(S_1;x_1,\ldots,x_{k+1})}\\
cdes_{n-1}(S_2) &=&\displaystyle{\sum_{\substack{ x_1, \cdots,
x_{k+1} \in
\{0,1\}}} R_{k+1}(S_2;x_1,\ldots,x_{k+1})}\\
cdes_{n-1}(S_3) &=&\displaystyle{\sum_{\substack{ x_1, \cdots, x_{k}
\in \{0,1\}}} R_{k}(S_3;x_1,\ldots,x_{k})}
\end{eqnarray*}

It is easy to check that
\begin{eqnarray*}&&\displaystyle{\sum_{\substack{ x_1, \cdots, x_{k+1} \in
\{0,1\}}} R_{k+1}(S;x_1,\ldots,x_{k+1})}-cdes_{n}(S_1)-cdes_{n-1}(S_2)\\
&=&\displaystyle{\sum_{\substack{ x_1, \cdots, x_{k+1} \in
\{0,1\}}}\frac{x_{k+1}-1}{1+x_1+\ldots
+x_{k+1}}R_{k+1}(S;x_1,\ldots,x_{k+1})}\\
&=&\displaystyle{\sum_{\substack{ x_1, \cdots, x_{k} \in
\{0,1\}}}-\frac{1}{1+x_1+\ldots
+x_{k}}R_{k+1}(S;x_1,\ldots,x_{k},0)}\\
&=&\displaystyle{\sum_{\substack{ x_1, \cdots, x_{k} \in
\{0,1\}}}R_{k}(S_3;x_1,\ldots,x_{k})}\\
&=&cdes_{n-1}(S_3).
\end{eqnarray*}

So, $\displaystyle{\sum_{\substack{ x_1, \cdots, x_{k+1} \in
\{0,1\}}}
R_{k+1}(S;x_1,\ldots,x_{k+1})}=cdes_{n}(S_1)+cdes_{n-1}(S_2)+cdes_{n-1}(S_3)$.
By Lemma \ref{DesR}, we have $ cdes_n(S) = cdes_n(S_1) +
cdes_{n-1}(S_2) + cdes_{n-1}(S_3)$. Hence,
$$cdes_n(S)=\displaystyle{\sum_{\substack{ x_1, \cdots, x_{k+1} \in
\{0,1\}}} R_{k+1}(S;x_1,\ldots,x_{k+1})}.$$ This complete the proof.
\end{proof}

We can associate $CDES_n(S)$ to a generating binary tree. We
consider the  generating binary tree as follows. $$\left \{
\begin{array}{ll}
                       \mbox{root} : & (1) \\
                       \mbox{rule} : & (k) \rightarrow (k)(k+1) \\
                      \end{array}
              \right.$$
This means first that the root has label $(1)$, and then for all k,
any node labeled $(k)$ will have $2$ descendants and they will have
the labels $(k),(k+1)$. We denote the obtained tree by $T$. For any
vertex $v$, there exists a unique path $P_{v}$ connecting $v$ and
the root. Define the height of $v$ to be the number of the edges in
$P_v$. If $\{w,v\}$ is an edge of $T$ and the height of a vertex $w$
is less than the height of $v$, then we say $w$ is the father of
$f$(resp. $v$ is a child of $w$), denoted by $w=f(v)$(resp. $v \in
c(w)$).

Suppose $E(T)$ is edge set of $T$. Define the edge weight function
of $T$, denoted $ew$, as $ew(\{w,v\})=(-1)^{1-k+l}$, where $(k)$ and
$(l)$ are the labels of the vertices $v$ and $w$ respectively, for
any $\{w,v\} \in E(T)$ with $w=f(v)$. Clearly,  $k-l\in\{0,1\}$.
Given a infinite sequence $\vec{d}=(d_1,d_2,\ldots,d_j,\ldots)$, we
define the vertex weight function of $T$, denoted $vw_{\vec{d}}$, as
$vw_{\vec{d}}(v)=k^{d_j}$, where $(k)$ is the labels of the vertex
$v$ and $j$ is the height of $v$, for any non-root vertex $v $, and
$vw(v)=1$ if $v$ is the root. For any $v\in V(T)$, recall that there
exists a unique path $P_{v}$ connecting $v$ and the root. We define
the weight $\omega_{\vec{d}}(P_v)$ of $P_v$ as
$\omega_{\vec{d}}(P_v)=\prod\limits_{u\in
V(P_v)}vw_{\vec{d}}(u)\prod\limits_{e\in E(P_v)}ew(e)$. For a fixed
$k\geq 1$, let $T_k$ be the subgraph of $T$ induced by the vertices
with height less or equal to $k$ and $Leaf(T_k)$ the set of the
leaves of $T_k$. Define the wight $\omega_{\vec{d}}(T_k)$ of $T_k$
as follows:
\begin{eqnarray*}\omega_{\vec{d}}(T_k)=\sum\limits_{v\in
Leaf(T_k)}\omega_{\vec{d}}(P_v)=\sum\limits_{ v \in
Leaf(T_k)}\prod\limits_{u\in
V(P_v)}vw_{\vec{d}}(u)\prod\limits_{e\in
E(P_v)}ew(e).\end{eqnarray*}

\begin{thm} Let $k\geq 1$ and $\vec{d}=(d_1,\ldots,d_k)$
be a sequence of nonnegative integers. Then
$$\omega_{\vec{d}}(T_k)=\sum\limits_{x_1,x_2,\ldots,x_k\in
\{0,1\}}(-1)^{k-\sum\limits_{i=1}^kx_i}\prod\limits_{i=0}^{k}\left(1+\sum\limits_{j=1}^ix_j\right)^{d_i}$$
\end{thm}
\begin{proof} Use $D$ to denote the set of all the $k$-tuples $(x_1,\ldots,x_k)$ such that $x_i\in \{0,1\}$.
 First, we establish
a bijection $\theta$ from $Leaf(T_k)$ to $D$. For any $v\in
Leaf(T_k)$, let $P_v$ be the unique path connecting $v$ and the
root. Suppose that $P_v=w_0w_1\ldots w_k$. Then $w_0$ is the root
and $w_k=v$. Let $x_i$ be the difference of the labels of $w_{i}$
and $w_{i-1}$ for any $i\in[k]$ and $\theta(v)=(x_1,\ldots,x_k)$. It
is easy to obtain $x_i \in [0,1]$ for any $i\in [k]$ by the
construction of the generating binary tree $T$. Therefore,
$\theta(v)\in D$. Conversely, for any $(x_1,\ldots,x_k)\in D$, let
$w_0=1$ and the label of $w_i$ is the sum the label of $w_{i-1}$ and
$x_i$ for any $i\in [k]$. Equivalently, the label of $w_i$ is
$1+\sum\limits_{j=1}^ix_i$ for all $i\in [0,k]$. So,
$\theta^{-1}((x_1,\ldots,x_k))=w_k \in Leaf(T_k)$. Hence, $\theta$
is a bijection from $Leaf(T_k)$ to $D$. Moreover, for any $v\in
Leaf(T_k)$, we have
\begin{eqnarray*}\omega_{\vec{d}}(P_v)&=&\prod\limits_{u\in V(P_v)}vw_{\vec{d}}(u)\prod\limits_{e\in E(P_v)}ew(e)\\
&=&\prod\limits_{i=1}^{k}\left(1+\sum\limits_{j=1}^ix_j\right)^{d_i}\prod\limits_{i=1}^{k}(-1)^{1-x_i}\\
&=&(-1)^{k-\sum\limits_{i=1}^kx_i}\prod\limits_{i=1}^{k}\left(1+\sum\limits_{j=1}^ix_j\right)^{d_i},\end{eqnarray*}
where $(x_1,\ldots,x_k)=\theta(v)$. Hence,
$$\omega_{\vec{d}}(\mathscr{T}_k)=\sum\limits_{x_1,x_2,\ldots,x_k\in
[0,1]}(-1)^{k-\sum\limits_{i=1}^kx_i}\prod\limits_{i=1}^{k}\left(r+\sum\limits_{j=1}^ix_j\right)^{d_i}.$$
\end{proof}

\begin{cor} Suppose that $n\geq 2$ and $S\subseteq [2,n]$. Let $\vec{d}=(d_1(S),\ldots,d_{|S|}(S))$. Then $\omega_{\vec{d}}(T_{|S|})=cdes_n(S).$
\end{cor}

For each $k \leq n$, consider the polynomial $\displaystyle{g_{n,k}
(\mathbf{x}) = \sum_{\sigma \in \mathfrak{S}_{n}, |CDES_n(\sigma)|=k
} \mathbf{x}_{CDES_n(\sigma)}}$. We call it {\em circular descent
polynomial with size $k$}. It is clear that $\displaystyle{ g_n =
\sum_{i=0}^k g_{n,k} ~  y^k }$. Let $\mathscr{D}_{k,n} = \{
(d_1,\cdots,d_k) :  d_1 + \cdots + d_k = n-1 ,d_i > 0 \mbox{ for all
i } \in [k] \}$. We define $\tau$ is a mapping from
$\mathscr{D}_{k,n}$ to $[2,n]$ as $\tau(D)=\{ 1+d_1, \cdots,
1+d_1+\cdots+d_k \}$ for each $D \in \mathscr{D}_{k,n}$. It is clear
that $\tau$ is the inverse mapping of $d$ in Theorem \ref{main}.
Similarly, we have the following corollary.

\begin{cor}
$\displaystyle{g_{n,k} (\mathbf{x}) = \sum_{D \in
\mathscr{D}_{n,k} } \omega( T_k) \mathbf{x}_{\tau(D)}}$ and
$\displaystyle{g_{n} (\mathbf{x},y) = \sum_{k=0}^n ( \sum_{D \in
\mathscr{D}_{n,k} } \omega( T_k) \mathbf{x}_{\tau(D)} ~ ) y^k }.$
\end{cor}

In the final of this section, we will give a example to show the
relation between $cdes_n(S)$ and $\omega_{\vec{d}}( T_k))$. Suppose
$k=2$, $n = 1+d_1+d_2$, $d=\{d_1,d_2\}$ and $S = \{ 1+d_1, 1+d_1+d_2
\}$. The weight $\omega_d(T_2)$ with respect to $d$ is Figure
\ref{T2}. Therefore, $\omega_{\vec{d}}(T_2)=
3^{d_2}2^{d_1}-2^{d_2+d_1}-2^{d_2}1^{d_1}+1^{d_1+d_2}$. By Theorem
\ref{main}, we have $\omega_{\vec{d}}(T_2)=cdes_n(S)$.
\begin{center}\label{T2}
\includegraphics[width=1.5in]{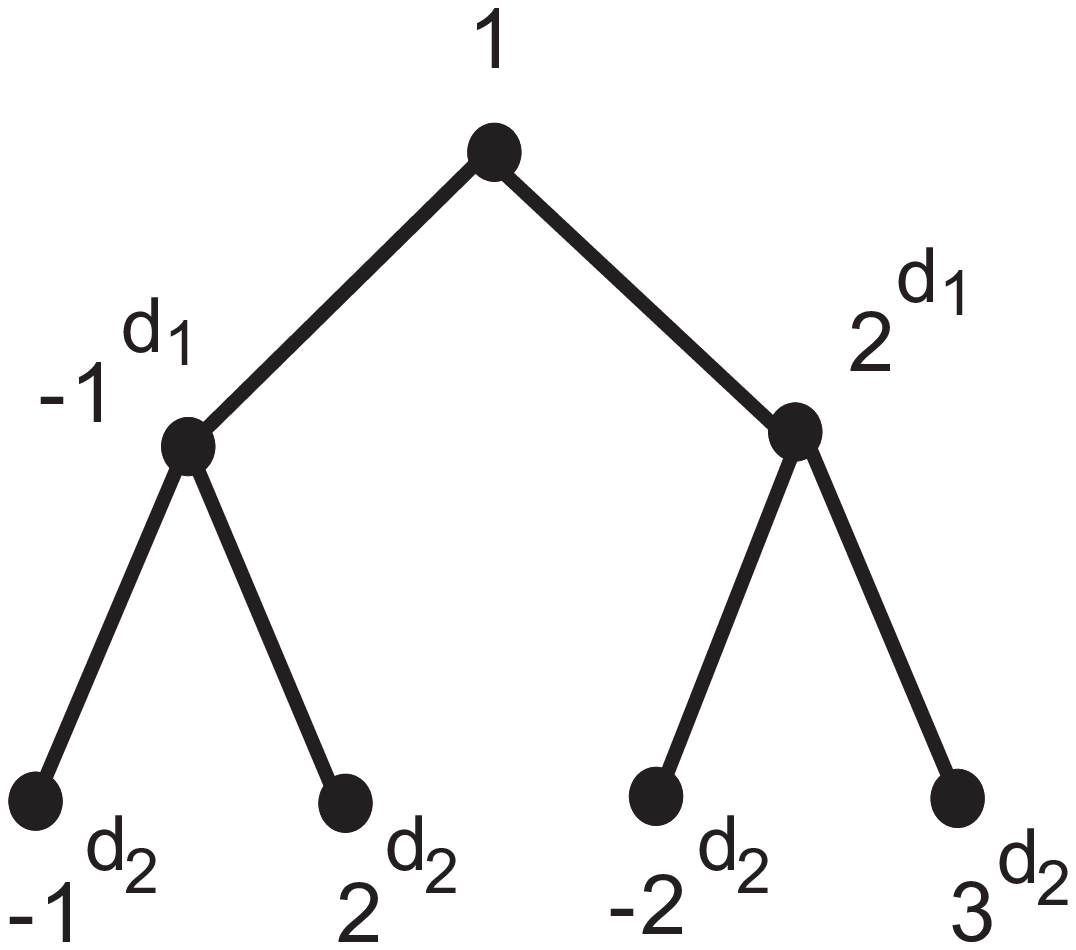}\\
Fig.\ref{T2}. The weight of the generating binary tree $T_2$
\end{center}

\section{Some Applications of The Main Theorem}
In this section, we will counting the number of permutation tableaux
according to their shape and give the enumeration of permutations by
the non-weak excedance bottom set.

 For any $S\subseteq [2,n]$, define the
{\it type} of the set $S$, denoted $type(S)$, as
$(r_1^{m_1},r_2^{m_2},\ldots,r_s^{m_s})$ if
$S=\bigcup\limits_{i=1}^{s}[r_i-m_i+1,r_i]$ such that $r_{i}\leq
r_{i-1}-m_{i-1}-1$ for all $i\geq 2$. Obviously,
$\sum\limits_{i=1}^{s}m_i=|S|$, $r_1=\max S$ and $r_s-m_s+1=\min S$.
We restate Theorem \ref{main} as follows.

\begin{lem}\label{lemmatype} Let $n$ be a positive integer with $n\geq 2$ and $S\subseteq [2,n]$.
Suppose that $type(S)=(r_1^{m_1},r_2^{m_2},\ldots,r_s^{m_s})$ and
let $r_{s+1}=1$ and $M_i=m_1+\ldots +m_i$ for any $i\in [s]$. Then
 \begin{eqnarray*}cdes_n(S) =
\displaystyle{\sum_{\substack{ x_1, \cdots, x_{|S|} \in \{0,1\}}}
(-1)^{|S|-\sum\limits_{j=1}^{|S|}x_j}\prod_{i=1}^{|S|}\left(1+\sum\limits_{j=1}^{i}x_j\right)
\prod_{i=1}^{s}
\left(1+\sum\limits_{j=1}^{M_i}x_j\right)^{r_{i}-r_{i+1}-m_{i}}
}.\end{eqnarray*}
\end{lem}

Suppose that the partition $Y_\lambda$ is contained in the
$k\times(n-k)$ rectangle. Regard the south-east border of
$Y_\lambda$ as giving a path $\mathcal {P}=\{P_i\}^n_{i=1}$ of
length $n$ from the northeast corner of the rectangle to the
southwest corner the rectangle: label each of the steps in this path
with a number from $1$ to $n$ according to the order in which the
step was taken. Sylvie Corteel and Philippe Nadeau \cite{SCPN}
established Bijection $I$ between permutation tableaux of length $n$
and permutations of $[n]$. We state as the following lemma.

\begin{lem}\label{lemmacdestab}{\rm \cite{SCPN}} There is a bijection $\Phi$ from permutation tableaux of length $n$ to $\mathfrak{S}_n$.
Moreover, for any permutation tableau $T_n^k$, the circular descent
of $\sigma=\Phi(T_n^k)$ are precisely the labels on the horizontal
edges of $\mathcal{P}$.
\end{lem}

For any a partition $\lambda=(\lambda_1,\ldots,\lambda_k)$, there
exists an unique strictly decreasing sequence
$a_\lambda=(a_1,\ldots,a_s)$ with $s\leq k$ such that there are some
$j$ satisfying $a_j=\lambda_i$ for any $i\in [k]$. On the other
hand, let $b_\lambda=(b_1,\ldots,b_s)$ such that $b_i=|\{j\in
[k]\mid \lambda_j\geq a_i\}|$ for all $i\in[s]$. The pair
$(a_\lambda,b_\lambda)$ is called the {\it type} of the partition
$\lambda$, denoted by $type(\lambda)$.

\begin{cor}Fix a partition $\lambda=(\lambda_1,\ldots,\lambda_k)$ with $\lambda_k\geq 1$.
Suppose that $type(\lambda)=(a_\lambda,b_\lambda)$, $b_0=1$ and the
length of $a_\lambda$ is $s$. Then the number of permutation
tableaux $\mathcal{T}^k_n$ with shape $\lambda$ is equal to
\begin{eqnarray*} \displaystyle{\sum_{\substack{ x_1, \cdots,
x_{n-k} \in \{0,1\}}}
(-1)^{n-k-\sum\limits_{j=1}^{n-k}x_j}\prod_{i=1}^{n-k}\left(1+\sum\limits_{j=1}^{i}x_j\right)
\prod_{i=1}^{s}
\left(1+\sum\limits_{j=1}^{a_{s+1-i}}x_j\right)^{b_{s+1-i}-b_{s-i}}
}.\end{eqnarray*}
\end{cor}
\begin{proof}
Let $\mathcal{T}_n^k$ be a permutation tableau with shape $\lambda$.
Since the type of $\lambda$ is $(a_\lambda,b_\lambda)$, let
$a_{s+1}=0$, all the labels of the horizontal edges of $\mathcal{P}$
are in the set $S=
\bigcup\limits_{i=1}^{s}[a_{1}-a_i+b_i+1,a_{1}-a_{i+1}+b_i]$.
Clearly, $|S|=\lambda_1=a_1=n-k$. Let $m_i= a_{s+1-i}-a_{s+2-i}$ and
$r_i=a_{1}-a_{s+2-i}+b_{s+1-i}$ for all $1\leq i\leq s.$ It is easy
to check that $type(S)=(r_1^{m_1},r_2^{m_2},\ldots,r_s^{m_s})$. From
Lemma \ref{lemmacdestab} , it follows that the number of permutation
tableaux $T^k_n$ with shape $\lambda$ is equal to $cdes_n(S)$. Lemma
\ref{lemmatype} implies that  the number of permutation tableaux
$T^k_n$ with shape $\lambda$ is
\begin{eqnarray*}
\displaystyle{\sum_{\substack{ x_1, \cdots, x_{n-k} \in \{0,1\}}}
(-1)^{n-k-\sum\limits_{j=1}^{n-k}x_j}\prod_{i=1}^{n-k}\left(1+\sum\limits_{j=1}^{i}x_j\right)
\prod_{i=1}^{s}
\left(1+\sum\limits_{j=1}^{a_{s+1-i}}x_j\right)^{b_{s+1-i}-b_{s-i}}
}.\end{eqnarray*}
\end{proof}

We also need the following lemma \cite{ESLK}.

\begin{lem}\label{lemmanonwektab}{\rm \cite{ESLK}} There is a bijection $\Psi$ from permutation tableaux of length $n$ to $\mathfrak{S}_n$.
Moreover, for any permutation tableau ${\cal T}_n^k$, the non-weak
excedance bottoms  of $\sigma=\Psi({\cal T}_n^k)$ are precisely
the labels on the horizontal edges of $\mathcal{P}$.
\end{lem}

Combining Lemmas \ref{lemmacdestab} and \ref{lemmanonwektab}, we
immediately obtain the following corollary.

\begin{cor} Let $n$ be a positive integer with $n\geq 2$ and $S\subseteq [2,n]$.
Then $nwexb_n(S)=cdes_n(S)$ and $\displaystyle{\sum_{\sigma \in
\mathfrak{S}_{n}} \mathbf{x}_{NWEXB_n(\sigma)} } y^{|NWEXB_n(\sigma)|} =\sum_{\sigma \in
\mathfrak{S}_{n}} \mathbf{x}_{CDES_n(\sigma)} y^{|CDES_n(\sigma)|} $.
\end{cor}



\end{document}